\documentclass[12pt]{amsart}
\usepackage{amsmath, amsthm, amsfonts}
\usepackage{amssymb}
\usepackage{amsbsy}
\usepackage{mathrsfs}

\input amssym.def
\input amssym
\usepackage{enumerate}
\usepackage{hyperref}
\hypersetup{
    pdftitle={Memory Estimation of Inverse Operators},    % title
    pdfauthor={Anatoly Baskakov}{Ilya Krishtal},     % author
    pdfsubject={Memory Estimation of Inverse Operators},   % subject of the document
    pdfcreator={Ilya Krishtal},   % creator of the document
    pdfproducer={Ilya Krishtal}, % producer of the document
    pdfkeywords={Banach modules}{Beurling spectrum}{Wiener's lemma}, % list of keywords
     colorlinks,%
    citecolor=blue,%
    filecolor=black,%
    linkcolor=magenta,%
    urlcolor=black
}

\chardef\bslash=`\\ % p. 424, TeXbook
\makeatletter
\def\verbatim{\interlinepenalty\@M \@verbatim
   \leftskip\@totalleftmargin\advance\leftskip2pc
   \frenchspacing\@vobeyspaces \@xverbatim}
\makeatother
\newtheorem{thm}{Theorem}[section]
\newtheorem{cor}[thm]{Corollary}
\newtheorem{lem}[thm]{Lemma}
\newtheorem{prop}[thm]{Proposition}

\theoremstyle{definition}
\newtheorem{defn}{Definition}[section]

\theoremstyle{remark}
\newtheorem{rem}{Remark}[section]
\newtheorem{exmp}{Example}[section]

\numberwithin{equation}{section}

\newcommand{\begeq}{\begin {equation}}
\newcommand{\eq}{\end{equation}}
\newcommand{\bs}{\begin {split}}
\newcommand{\es}{\end{split}}
\newcommand{\bp}{\begin {prop}}
\newcommand{\ep}{\end {prop}}
\newcommand{\bt}{\begin {thm}}
\newcommand{\et}{\end {thm}}
\newcommand{\bc}{\begin {cor}}
\newcommand{\ec}{\end {cor}}
\newcommand{\bl}{\begin {lem}}
\newcommand{\el}{\end {lem}}
\newcommand{\bpf}{\begin {proof}}
\newcommand{\epf}{\end {proof}}
\newcommand{\bi}{\begin {itemize}}
\newcommand{\ei}{\end {itemize}}
\newcommand{\ben}{\begin {enumerate}}
\newcommand{\een}{\end {enumerate}}
\newcommand{\brem}{\begin {rem}}
\newcommand{\erem}{\end {rem}}

\newcommand{\bd}{\begin {defn}}
\newcommand{\ed}{\end {defn}}

\newcommand{\bex}{\begin {exmp}}
\newcommand{\eex}{\end {exmp}}

\newcommand{\B}{{\mathcal B}}

\newcommand{\F}{\mathcal{F}}

\newcommand{\TTT}{{\T\kern-.44em \T}}
\newcommand{\tTTT}{\widetilde{\T\kern-.44em \T}}

\newcommand{\ZZ}{{\mathbb Z}}
\newcommand{\TT}{{\mathbb T}}
\newcommand{\RR}{{\mathbb R}}
\newcommand{\CC}{{\mathbb C}}
\newcommand{\NN}{{\mathbb N}}

\newcommand{\X}{\mathcal{X}}

\newcommand{\s}{\sigma}
\renewcommand{\l}{\lambda}
\renewcommand{\L}{\Lambda}
\renewcommand{\a}{\alpha}

\renewcommand{\SS}{\mathcal S}

\newcommand{\T}{\mathcal{T}}

\newcommand{\LL}{\mathscr{L}}

\DeclareMathOperator{\Imm}{Im }

\DeclareMathOperator{\supp}{supp}
\newcommand{\esssup}{\text{ess \!sup }\!}

\begin{document}

\title[Parabolic Operators in Homogeneous Function Spaces]{Harmonic and Spectral Analysis of Abstract Parabolic Operators in Homogeneous Function Spaces} 
%{Gerhard-Pr\"uss Theorem in Homogeneous Banach Spaces}%{Frames, Fusion Frames, Fuzzy Fusion Frames and Wiener's Lemma}

%%%%%%%%%%%%%%%%%%%%%%%%%%%%%%%%%%%%%%%%%%%%%%%%%%%%%%%%%%%%%%%%%%%%%%%%
\author{Anatoly G. Baskakov and Ilya A. Krishtal}
\address{Department of Applied Mathematics and Mechanics, Voronezh State University, Russia 394693 \\
email: mmio@amm.vsu.ru}
\address{Department of Mathematical Sciences, Northern Illinois University, DeKalb, IL 60115 \\
email: krishtal@math.niu.edu}

\thanks{ The first author is supported in part by RFBR grant 13-01-00378.
The second author is supported in part by NSF grant DMS-1322127.}

%%%%%%%%%%%%%%%%%%%%%%%%%%%%%%%%%%%%%%%%%%%%%%%%%%%%%%%%%%%%%%%%%%%%%%%%

\date{\today }

%\subjclass[2010]{46H25}

\keywords{Abstract parabolic operators, Homogeneous function spaces, Green functions, Beurling spectrum}

%\dedicatory{to }

%\thanks{Research of the first author was supported in part by NSF
%{grant DMS-9805483.}
%%%%%%%%%%%%%%%%%%%%%%%%%%%%%%%%%%%%%%%%%%%%%%%%%%%%%%%%%%%%%%%%%%%%%%%%
%%%%%%%%%%%%%%%%%%%%%%%%%%%%%%%%%%%%%%%
%%%%%%%%%%%%%%%%%%%%%%%%%%%%%%%%%
\begin{abstract}
We use methods of harmonic analysis and group representation theory to study the spectral properties
of the abstract parabolic operator $\LL = -d/dt+A$ in homogeneous function spaces. We provide 
 sufficient conditions for invertibility of such operators in terms of the spectral properties of the operator $A$ and the semigroup generated by $A$. We introduce a homogeneous space of functions with absolutely summable
 spectrum and prove a generalization of the Gearhart-Pr\"uss Theorem for such spaces. We use the results
 to prove existence and uniqueness of solutions of a certain class of non-linear equations.
\end{abstract}
%%%%%%%%%%%%%%%%%%%%%%%%%%%%%%%%%%%%%%%
%%%%%%%%%%%%%%%%%%%%%%%%%%%%%%%%%
\maketitle

\section{Introduction}
%\baselinestretch{1.5}
In this paper we study the spectral properties of a differential operator %(see Definition \eqref{defdifop})
\begeq\label{difop}
\LL = -d/dt + A:\ D(\LL)\subset \F(\RR, X) \to \F(\RR, X)
\eq
in homogeneous Banach spaces $\F(\RR, X)$ %(see Definition \ref{defhom}) 
of functions with values in a complex Banach space $X$. The operator $A: D(A)\subset X\to X$ in \eqref{difop} is assumed to be the infinitesimal generator of a $C_0$-semigroup $T: \RR_+ = [0,\infty) \to B(X)$. The homogeneous spaces $\F =\F(\RR,X)$
and the operator $\LL$ are identified precisely in Definitions \ref{defhom} and \ref{defdifop}, respectively.

The properties of the homogeneous spaces $\F$ allow us to correctly define the \emph{Howland semigroup}
$\mathscr T: \RR^+ \to \B(\F)$ by
\[(\mathscr T(t)x)(s) = T(t)x(s-t), \  s\in\RR, x\in\F, t\in\RR^+.\]

In some homogeneous spaces such as $C_0(\RR, X)$ and $L^p(\RR, X)$, see Example \ref{exhom} for the definitions, it was proved \cite{B95, B96} that the operator $\LL$ is the infinitesimal generator of the semigroup
$\mathscr T$. Moreover, in a large class of homogeneous spaces the following result holds.

\bt \label{thm1}
The following are equivalent:
\begin{itemize}
\item The operator $\LL$ is invertible, that is, the equation
\begeq\label{difeq1}
\frac {dx}{dt} = Ax+y
\eq
has a unique (mild) solution $x\in\F(\RR, X)$ for any $y\in \F(\RR, X)$;
\item The semigroup $T$ is \emph{hyperbolic}, that is, the spectrum $\s(T)$ of the semigroup $T$ satisfies $\s(T(1))\cap \TT = \emptyset$, $\TT = \{\l\in\CC:\
|\l| = 1\}$;
\item The Howland semigroup $\mathscr T$ is hyperbolic.
\end{itemize}
\et

For $\F\in \{ C_0, L^p\}$, $1\le p <\infty$,  the above result was proved in \cite{LMS94, LMS95} and \cite[Theorem 2.39]{CL99}. For a larger class of spaces, including $L^\infty$ and $C_b$, the theorem appears
in \cite{B95, B96, B13}. We remark that in many homogeneous spaces the result still holds even though the Howland semigroup $\T$ is not strongly continuous (see \cite{B13}).

In this paper we prove that one of the implications in the above theorem holds for all homogeneous spaces in Definition \ref{defhom}:

\bt
If $T$ is a hyperbolic semigroup,  then the operator $\LL$  is invertible. 
\et

If $X$ is a Hilbert space, classical results of L.~Gearhart and J.~Pr\"uss \cite{Ge78, P84, EN00, P01} provide another equivalence to the statements in Theorem \ref{thm1}.

\bt
The operator $\LL$ is invertible if and only if the imaginary axis $i\RR$ is a subset of the resolvent set $\rho(A)$ of the generator $A$, i.e.
\begeq\label{intera}
\s(A) \cap (i\RR) = \emptyset,
\eq 
 and the resolvent operator
$R(\l, A) = (A -\l I)^{-1}$ satisfies 
\begeq\label{resbou}
M  = \sup\limits_{\l\in i\RR} \|R(\l, A)\| <\infty.
\eq
\et

If $X$ is not a Hilbert space, this  equivalence does not hold in general (see \cite{GVW81} and \cite[Counterexample IV.2.7]{EN00}). Typical generalizations of the above result to Banach spaces \cite[and references therein]{LS01} would impose additional restrictions on the resolvent $R(\l, A)$ or the Banach space $X$. 
In this paper, 
 we pursue a different kind of generalization, where we deal with an arbitrary (complex) Banach space $X$
and impose no additional restrictions on the resolvent. Instead, 
we define a class $\F_{as}$ of functions in $\F(\RR, X)$ which is a Banach function space where the differential operator $\LL$ is invertible provided that \eqref{intera} and \eqref{resbou} hold. The space $\F_{as}$ of functions with the absolutely summable spectrum (see Definition \ref{asdef}) is defined using the spectral theory of Banach modules \cite{B04, HR79}.

We refer to \cite{ABHN11, CL99, EN00} for more information on the background and history of research related
to this paper.

The remainder of the paper is organized as follows. In Section \ref{sec:prelim} we define homogeneous function spaces and  introduce the space $\F_{as}$ of functions with the absolutely summable spectrum. We do the latter by means
of the spectral theory of Banach $L^1(\RR)$-modules, basic notions of which are also discussed in  Section
\ref{sec:prelim}. In Section \ref{operat} we study the basic properties of the differential operator $\LL_\F$ and
prove the first of our two main results: a sufficient condition for invertibility of
$\LL_\F$ in a homogeneous Banach function space, which appears in Theorem \ref{gfthm}. Our other main result, Theorem \ref{mainthm}, appears in Section \ref{sec:4} and establishes sufficient conditions for invertibility of the
operator $\LL$ in any homogeneous space of functions with the absolutely summable spectrum. We also provide an estimate for $\|\LL_{\F_{as}}^{-1}\|$. In the final section of the paper we illustrate our results with a  counterexample to the
Gearhart-Pr\"uss Theorem and an application  to a special kind of non-linear differential equations.

\section{Preliminaries}\label{sec:prelim}
In this section we introduce the notation, define homogeneous function spaces, and survey the necessary tools from the spectral theory of Banach $L^1(\RR)$-modules.  We also introduce the space of functions with the absolutely summable spectrum.

The symbol $X$ %and $Y$ 
will denote a complex Banach space and $B(X)$ will be the Banach algebra of
all bounded linear operators on $X$. By $T: \RR_+\to B(X)$ we shall denote a $C_0$-semigroup of operators in $B(X)$ and $A: D(A)\subset X\to X$ will be its infinitesimal generator \cite{EN00}.

By $L^1(\RR,X)$ we shall denote the Banach space of all (equivalence classes) of Bochner integrable $X$-valued functions with the standard $L^1$-norm:
\[\|f\| = \|f\|_1 = \int_\RR \|f(t)\|_{X} dt,\quad f\in L^1(\RR,X).\]

If $X = \CC$, we shall use the notation $L^1 = L^1(\RR)$ for the standard group algebra
of %(equivalence classes of) 
Lebesgue integrable functions. For $f \in L^1(\RR)$, we shall denote by $\hat f$ the Fourier transform of $f$ given by
\[
\hat f(\l) = \int_\RR f(t)e^{-it\l} dt, \ \l\in\RR.
\]

The space $L^1(\RR, B(X))$ defined as above, however, may be too small for our purposes. Occasionally, we will use the space $L^1_s(\RR, B(X))$ that consists
of all %(equivalence classes of) 
%of Bochner integrable $B(\X)$-valued 
functions $F: \RR \to B(X)$ with the following properties:
\ben
\item For all $x\in X$ the function $s\mapsto F(s)x: \RR\to X$ is measurable;
\item There is $f\in L^1(\RR)$ such that 
\begeq\label{normmajor}\|F(s)\|\le f(s)\ \mbox{a.e.}\eq %for a.e. $s\in\RR$.
\een
For $F\in L^1_s(\RR, B(X))$ we let
\[\|F\| = \|F\|_1 = \inf \|f\|,\]
where the infimum is taken over all functions $f$ that satisfy \eqref{normmajor}.

The space $L^1_s(\RR, B(X))$ is a Banach algebra with the
multiplication %in $L^1(\RR, B(\X))$ is played by the convolution of functions  
given by
\[(F_1*F_2)(t) = \int_\RR F_1(s)F_2(t-s)ds, \quad F_1, F_2\in L^1_s(\RR, B(X)).\] 
In particular, $\|F_1*F_2\| \le \|F\|_1\|F\|_2$.

We shall also use the space of locally integrable $X$-valued functions $L^1_{loc}(\RR,X)$, which consists of all  measurable functions $f: \RR \to X$ such that
\[\int_K \|f(t)\|_{X}dt < \infty\]
for any compact set $K\subset \RR$.

For $p\in[1,\infty)$, the Stepanov space $\SS^p = \SS^p(\RR,X)$ consists of all functions  $x\in L^1_{loc}(\RR,X)$ such that
\[\|x\|_{\SS^p} = \sup_{t\in\RR} \left(\int_0^1 \|x(s+t)\|^p ds\right)^{1/p} < \infty.\]

%\subsection{
%Banach $\A(\X)$-modules}\label{banmod}

\subsection{Homogeneous Banach function spaces
%$\A(\X)$-modules
}\label{homspac}\

In this paper we consider  Banach function spaces $\F(\RR, X)$ that are homogeneous according to the following definition.

\bd\label{defhom}
A Banach function space $\F = \F(\RR, X)$ is \emph{homogeneous} if it has the following properties:
\ben
\item %$\F(\X) \subseteq L^1_{loc}(\RR,\X)$.
%We have $\F(\RR, \X) \subseteq S^1(\RR,\X)$, that is 
$\F$ is continuously embedded into $\mathcal S^1$;
\item For all $t\in\RR$ and $x\in \F$ we have $S(t)x \in \F$, where
\begeq\label{transrep}
S(t)x(s) = x(t+s), \eq
and the translation operator $S(t)$ is an isometry in $B(\F)$.
\item Given $x\in \F$ and $C \in  B(X)$, the  function
\[ y(t) = C(x(t))\]
belongs to $\F$ and $\|y\| \le\|C\| \|x\|$.
%The \emph{translation} representation $T: \RR \to B(\F(\X))$ given by
%\begeq\label{trans1}
%T(t)f(s) = f(s+t),\ f\in\F, \ s,t\in\G,
%\eq
%is bounded. 
\item\label{svert} Given $x\in \F$ and $F \in L^1_s(\RR, B(X))$, the convolution
\[(F*x)(t) = \int_\RR F(s)x(t-s)ds\] %,\quad A\in \A(\X),\ f\in\F(\X),\]
belongs to $\F$ and $\|F * x\| \le\|F\|_1\|x\|$.
%defines a Banach $\A(\X)$-module structure on $\F(X)$.
\item If $x \in \F$ is such that $f*x = 0$ for all $f \in L^1(\RR)$ then $x = 0$.
\een
\ed

%Let us consider several examples of homogeneous Banach spaces. %Verification of the

\bex\label{exhom}
The following Banach spaces are homogeneous or have an equivalent norm that makes them homogeneous:
\ben
\item The spaces $L^p = L^p(\RR, X)$, $p\in [1,\infty]$, of functions $x\in L^1_{loc}(\RR,X)$ such that
\[\|x\|_{L^p} = \left(\int_\RR \|x(s)\|^p ds\right)^{1/p} < \infty, \ p\in [1,\infty),\]
or $\|x\|_\infty = \underset{t\in\RR}\esssup |x(t)| < \infty;$
\item Stepanov spaces $\SS^p = \SS^p(\RR,X)$, $p\in [1,\infty)$;
\item Wiener amalgam spaces $L^{p,q} = L^{p,q}(\RR, X)$, $p,q\in [1,\infty)$, of functions $x\in L^1_{loc}(\RR,X)$ such that
\[\|x\|_{L^{p,q}} = \left(\sum_{k\in \ZZ} \left(\int_0^1 \|x(s+k)\|^p ds\right)^{q/p}\right)^{1/q} < \infty, \ p,q\in [1,\infty);\]
\item The space $C_b = C_b(\RR, X)$ of bounded continuous $X$-valued functions with the norm
\[\|x\|_\infty = \sup_{t\in\RR} \|x(t)\|, \ x\in C_b;\]
\item The subspace $C_{ub} = C_{ub}(\RR, X) \subset C_b$ of uniformly continuous functions;
\item The subspace $C_{0} = C_{0}(\RR, X) \subset C_{ub}$ of continuous functions vanishing at infinity:  
$x\in C_0$ if $\lim_{|t|\to \infty} \|x(t)\| = 0$;
\item The subspace $C_{sl,\infty}= C_{sl,\infty}(\RR, X) \subset C_{ub}$ of slowly  varying at infinity functions: 
$x\in C_{sl,\infty}$ if $\lim_{|\tau|\to \infty} \|x(\tau+t) - x(\tau)\| = 0$ for all $t\in\RR$ (see \cite{BKal12});
\item The subspace $C_\omega = C_\omega(\RR, X) \subset C_{ub}$ of $\omega$-periodic functions, $\omega\in\RR$; 
\item The subspace $AP = AP(\RR, X) \subset C_{ub}$ of (Bohr) almost periodic functions \cite{BK10, LZ82}; 
\item The subspace $AP_{\infty} = AP_{\infty}(\RR, X) \subset C_{ub}$ of almost periodic at infinity functions \cite{B13} defined by
\[AP_{\infty} = \overline{\rm span}\{e^{i\l \cdot}x: \l\in\RR, x\in C_{sl,\infty}\}.\]
\item The spaces $C^k = C^k(\RR, X)$, $k\in \NN$, of $k$ times continuously differentiable functions with a bounded $k$-th derivative and the norm
\[
\|x\|_{(k)} = \|x\|_{\infty} + \|x^{(k)}\|_\infty < \infty;
\]
\item The H\"older spaces $C^{k,\a} = C^{k,\a}(\RR, X)$, $k\in\NN\cup\{0\}$, $\alpha\in (0,1]$:
\[
C^{k,\a} = \left\{x\in C^{k}:\  \|x^{(k)}\|_{C^{0,\a}} = \sup_{t\neq s\in\RR} \frac{|x(t)-x(s)|}{|t-s|^\a}<\infty\right\}, 
\]  
\[
\|x\|_{C^{k,\a}} = \|x\|_{C^k}+  \|x^{(k)}\|_{C^{0,\a}}.
\]
%\item Periodic at infinity?
%\item Modulation spaces?
\een
\eex

\brem
We note that  Definition \ref{defhom} in this paper differs from Definition 2.1 in \cite{B13}, which is more narrow. In particular, here we do not assume that the space $\F$ is solid and this allows us to consider the spaces of periodic and almost periodic functions.
\erem

\bd\label{sphom}
A homogeneous space $\F = \F(\RR, X)$ is called \emph{spectrally homogeneous} if for all $\l\in \RR$ and
$x\in \F$ we have $e^{i\l\cdot}x \in \F$ and $\|e^{i\l\cdot}x\| = \|x\|$.
\ed

Among all of the homogeneous spaces in Example \ref{exhom} the only class of spaces that are not spectrally homogeneous are the spaces of $C_\omega$ of $\omega$-periodic functions.

In a spectrally homogeneous space $\F$ there is a well-defined isometric representation $V\in \RR\to B(\F)$ given by 
\begeq\label{repV}
V(\l)x(t) = e^{i\l t}x(t),\ \l, t\in\RR, x\in\F.
\eq 
%$\l, t\in\RR$, $x\in\F$.

\subsection{
Banach $L^1$-modules and the Beurling spectrum}\label{banmod}\ 

Properties (4) and (5) in Definition \ref{defhom} ensure that a homogeneous Banach space is a non-degenerate
Banach $L^1(\RR)$-module \cite{B04, HR79}.
In this subsection we present the necessary definitions and results from the spectral theory of such Banach modules.
The  proofs  omitted in the presentation and further details can be found in
 \cite{B04, BK05, BK14, HR79}.
 
 \bd\label{defbanmod}
 A Banach space $\X$ is a \emph{Banach $L^1(\RR)$-module} if 
 there is a bilinear map $(f,x) \mapsto fx: L^1(\RR)\times \X\to \X$ such that
 \ben
 \item $(f*g)x = f(gx)$, $f,g \in L^1(\RR)$, $x\in\X$;
 \item $\|fx\| \le \|f\|\|x\|$, $f \in L^1(\RR)$, $x\in\X$.
 
\noindent The module structure is \emph{non-degenerate} if, in addition,
 \item $fx = 0$ for all $f \in L^1(\RR)$ implies $x = 0\in\X$.
  \een
 \ed

We say that the structure of a Banach $L^1(\RR)$-module $\X$ is  \emph{associated with a representation} 
$U: \RR\to B(\X)$ if
\begeq\label{assoc}
U(t)(fx) = f_t x = f(U(t)x), \ f \in L^1(\RR),\ x\in \X, \ f_t(s) = f(s+t).
\eq

As we mentioned above, any homogeneous function space $\X = \F(\RR, X)$ is a non-degenerate
Banach $L^1(\RR)$-module. Its structure is given by
\begeq\label{modstruc}
(fx)(t) = (f*x)(t) =\int_\RR f(t-s)x(s) ds = \int_\RR f(s)x(t-s) ds,
\eq
$f \in L^1(\RR)$, $x\in \F$, and is associated with the translation representation 
 $S$ defined by \eqref{transrep}.

\bd\label{bspec}
The \emph{Beurling spectrum} of an element $x \in \X$ is the subset $\L(x)\subseteq\RR$ the complement of which is given by
\[\{\l\in\RR: \mbox{ there is } f\in L^1(\RR) \mbox{ such that } \hat f(\l) \ne 0 %\]
%\[
\mbox{ and } fx =  0\}.\]
\ed

\brem
In homogeneous Banach spaces the Beurling spectrum $\L(x)$ coincides with the support of the (distributional) Fourier transform of $x\in\F$.
\erem

In the next lemma we present basic properties of the Beurling spectrum that will be  used throughout the paper. We refer to \cite{B04, BK05, R68} and references therein for the proof.

\bl\label{sprop}
Let $\X$ be a non-degenerate Banach $L^1(\RR)$-module. %with the structure associated with a representation $S$. 
Then 
\begin{description}
\item[ (i)]   $\L(x)$ is closed for every $x \in \X$ and $\L(x)=\emptyset$ if and only if $x =0$;
\item[(ii)] $\L(Ax + By) \subseteq \L(x) \cup \L(y)$ for all $A$, $B \in B(\X)$ that commute with all 
operators $x\mapsto fx$,
$f \in L^1(\RR)$; 
\item[(iii)] $\L(f x)\subseteq (\supp \hat f)\cap \L(x)$ for all $f \in L^1(\RR)$ and $x\in\X$;
%\item[(iv)] $f x=0$ if $(\supp \hat f)\cap\L(x)=\emptyset$, where $f \in L^1(\RR)$ and $x\in\X$;
%\item[(v)] $f x = x$ if $\L(x)$ is a compact set, and $\hat f \equiv 1$ in some neighborhood of $\L(x)$, $f\in L^1(\RR)$, $x\in\X$.
\item[(iv)] $f x=0$ if $(\supp \hat f)\cap\L(x)$ is countable and $\hat f(\l) = 0$ for all $\l \in (\supp \hat f)\cap\L(x)$, %a strong Ditkin set with empty interior, 
$f \in L^1(\RR)$, $x\in\X$;
\item[(v)] $f x = x$ if $\L(x)$ is a compact set, the boundary of $\L(x)$ is countable, and $\hat f \equiv 1$ on  $\L(x)$, $f\in L^1(\RR)$, $x\in\X$.
%\item[(vi)] if $M_0$ is dense in $M\subseteq X$, then $\L(M) = \overline{\bigcup_{x\in M_0}\L(x)}$.
\end{description}
\el

Given a closed set $\Delta\subset \RR$ we shall denote by $\X(\Delta)$ the (closed) \emph{spectral submodule} of
all vectors $x\in\X$ such that $\L(x)\subseteq \Delta$. The symbol $\X_{Comp}$ will stand for the set
of all vectors $x$ such that $\L(x)$ is compact. If the module structure is associated with a representation $U$, by $\X_U$ we shall denote the submodule of $U$-continuous vectors, i.e., the set of all vectors $x\in\X$ such that the function $t\mapsto U(t)x:\RR\to \X$ is continuous.

\brem
Observe that any spectral submodule $\F(\Delta)$ of a homogeneous function space $\F(\RR, X)$ is itself a homogeneous function space. It may happen that $\F(\Delta) = \{0\}$ even if $\Delta \neq \emptyset$. For example, if $\F = L^p(\RR, X)$, $1\le p <\infty$, and $\Delta$ is finite then $\F(\Delta) = \{0\}$. If $\Delta$ is compact then $\F(\Delta)\subset C_{ub}$ and each $x\in \F(\Delta)$ extends to an entire function of exponential type $\omega = \max\{|\l|, \l\in\Delta\}$ \cite{BKal12}.
\erem

The following lemma opens up the possibility of applying our main results to non-linear equations.

\bl
Let $\X$ be a non-degenerate Banach $L^1(\RR)$-module. Assume that the module structure is associated with a representation $U$ as in \eqref{assoc}. %\footnote{ Here and in the following subsection the group representation $S$ does not have to be the translation.}. 
Let $F: \X^n \to \X$, $\X^n =\underset{n \ times%\mbox{ \small{times}}
}{\underbrace{\X\times\dots\times\X}}$, be an $n$-linear map such that for any %$f\in L^1(\RR)$ and $k = 1, \ldots, n$ 
$t\in\RR$ we have
%\[f(F(x_1,\ldots, x_n)) = F(fx_1,\ldots, fx_k, \ldots, fx_n).\]
\[U(t)(F(x_1,\ldots, x_n)) = F(U(t)x_1,\ldots, U(t)x_n).\]
Then
\[\L(F(x_1,\ldots, x_n)) \subseteq \overline{\L(x_1)+\ldots +\L(x_n)}.\]
\el

\bpf
In case $\X$ is a Banach algebra and $F(x_1,x_2) = x_1x_2$, this result appears in \cite{BK05, BK14}. In the general case the proof may be repeated nearly verbatim with obvious modifications.
\epf

\subsection{%The space of 
Vectors with the absolutely summable spectrum}\

In this subsection we define the class $\X_{as}$ of vectors in a Banach $L^1(\RR)$-module $\X$ that have absolutely summable spectrum. In the scalar case this class was introduced in \cite{B97Sib}. In this exposition we follow \cite{BK14} where the general definition appears. We use the family of functions $(\phi_\a)$,
$\a\in\RR$, defined via the Fourier transform by
\begeq\label{triangled2}
\hat\phi_{a} (\l) = \hat\phi(\l-a),\  a\in\RR, %\ \mbox{ where }
%\hat \phi_{N}^d (\l_1,\ldots,\l_d)= \prod_{k=1}^d \hat\phi_{N}(\l_k),  
\eq
where %the functions $\phi_N$, $N\in\NN$, are defined via
\begeq\label{triangle}
%\hat\phi_N(\l) = \hat\phi(\frac\l N),\quad
\hat\phi(\l) \equiv \hat\phi_0(\l)=(1-|\l|)\chi_{[-1,1]}(\l),
\eq
and $\chi_E$ is, as usually, the characteristic function of the set $E$.

\bd\label{asdef} The class $\X_{as}$ of vectors in a Banach $L^1(\RR)$-module $\X$ that have \emph{absolutely summable spectrum} is 
\begeq\label{windef}
\X_{as}= \left\{x\in\X:\ \|x\|_{\widetilde{as}} =\int_{\RR^d}\|\phi_{a} x\|da< \infty \right\},
\eq
where the functions $\phi_{a}$ %arise from a standard family of BUPU via  
are defined by \eqref{triangled2}  and
\eqref{triangle}.
\ed

In \cite{BK14}, one can find plenty of examples of the spaces with absolutely summable spectrum. Classical Wiener amalgam spaces \cite{F80, FS85} are among the well-studied spaces that arise in such a way. 

\brem
In the above definition, instead of the family of functions $(\phi_{a})$ one can use just about any bounded uniform partition of unity \cite{F80, FG85}. One would obtain the same space as a result \cite{BK14}. This is analogous to using different window functions in the short time Fourier transform \cite{G01}.
\erem

In \cite{BK14} we have shown that 
$\X_{as}$ is a Banach space with the norm $\|\cdot\|_{\widetilde{as}}$. For our purposes, however, it is often more convenient to use an equivalent 
 norm given by
\begeq\label{wiensumel}
\|x\|_{as} = 5 \sum_{n\in\ZZ} \|\phi_{n} x\|.
\eq
In \cite{BK14} we obtained the inequalities  
%\begeq\label{normeq}
%\frac1{(2N+1)} \|x\|_{as,1} \le \|x\|_{as,N} \le 3\|x\|_{as,1},\ \mbox{and}
%\eq
\begeq\label{normeq1}
\|x\|_{\widetilde{as}} \le \|x\|_{as} \le 20\|x\|_{\widetilde{as}},\ x\in\X_{as}.
\eq
In \cite{BK14} we have also shown that if $\X$ is a Banach algebra and the module structure
is associated with a group of algebra automorphisms, then $\X_{as}$ is also a Banach algebra and, due
to the choice of the constant in \eqref{wiensumel},
\begeq\label{asalg}
\|xy\|_{as} \le\|x\|_{as}\|y\|_{as},\ x,y\in\X_{as}.
\eq
Moreover, the key result of \cite{BK14} states that the algebra $\X_{as}$ is inverse closed, i.e.~if $x\in\X_{as}$
is invertible in $\X$, then $x^{-1}\in \X_{as}$. The following analog of \eqref{asalg} can be proved in exactly the same way.

\bp
Let $\X$ be a non-degenerate Banach $L^1(\RR)$-module. Assume that the module structure is associated with a representation $U$ as in \eqref{assoc}. 
Let $F: \X^n \to \X$, be an $n$-linear map such that for any %$f\in L^1(\RR)$ and $k = 1, \ldots, n$ 
$t\in\RR$ we have
%\[f(F(x_1,\ldots, x_n)) = F(fx_1,\ldots, fx_k, \ldots, fx_n).\]
\[U(t)(F(x_1,\ldots, x_n)) = F(U(t)x_1,\ldots, U(t)x_n).\]
Then 
\begeq
\|F(x_1,\ldots, x_n)\|_{as} \le \|F\|\cdot \prod_{k=1}^n\|x_k\|_{as}.
\eq
\ep

It has been observed by many people, see, e.g., \cite{GK10} and \cite[Remark 3.5]{BK14}, that %if $\X = C_{ub}(\RR, X)$, 
smoothness of the function $x$ %(the property $x\in C^{1,\a}$) 
is closely related to the spectral decay of $x$. In particular, %for $x\in C_{ub}(\RR, H)$, where $H$ 
%if $X$ is a Hilbert space, 
we have   $C^{1,\a}\subset \X_{as}$, $\a>0$. 
Below we prove a slightly
weaker sufficient condition that uses the following modulus of continuity.

\bd Let $\X$ be a non-degenerate Banach $L^1(\RR)$-module with the structure associated with an isometric representation $U: \RR \to B(\X)$. For $x\in\X$, its \emph{modulus of continuity} $\omega_x$ is 
defined by
\[\omega_x (t)= \sup_{|s|\le t} \|U(s)x-x\|,\ t\ge 0.\]
\ed

\brem
The basic properties of the modulus of continuity can be found, for example in \cite{BB67, K76}. Here we mention 
the obvious facts that 
\[\lim_{t\to 0} \omega_x(t) = 0, x\in \X_U,\]
and $\omega_x$ is subadditive, that is
\[\omega_x(t+s)\le \omega_x(t)+\omega_x(s),\ x\in\X.\]
As an immediate consequence of subbaditivity we get
\begeq\label{modest}
\omega_x(ks)\le k\omega_x(s)\quad\mbox{and}\quad\frac{\omega_x(1)}k \le \omega_x(\frac1k),\ x\in\X, k\in\NN, s\ge 0.
\eq
Using the monotonicity  of $\omega_x$  we also get
\begeq\label{modestt}
\omega_x(\l s)\le (\l+1)\omega_x(s), \ \l,s> 0,\ x\in\X.
\eq
\erem 

\bl For $x\in\X$ we have
\begeq\label{modest1}
\|\phi_{k}x\| \le Const\cdot \omega_x(\frac1{|k|}), k\in\ZZ\backslash\{0\}.
\eq
\el

\bpf
Observe that for $k\neq 0$
\[\phi_{k}x = \int_\RR \phi(t)e^{ikt}U(-t)xdt = -\int_\RR \phi(\tau+\frac\pi k)e^{ikt}U(-\tau-\frac\pi k)xd\tau,\]
where we made the change of variables $t = \tau+\pi/k$. Averaging the above two expressions we get
\[\begin{split}
\phi_{k}x &= \frac12\int_\RR e^{ikt}(\phi(t)U(-t) - \phi(t+\frac\pi k)U(-t-\frac\pi k))xdt \\
&=\frac12\int_\RR e^{ikt}(\phi(t) - \phi(t+\frac\pi k))U(-t)xdt \\
&+\frac12\int_\RR e^{ikt}\phi(t+\frac\pi k)(U(-t) - U(-t-\frac\pi k))xdt.
\end{split}
\]
Hence,
\[
\|\phi_{k}x\| \le \frac12\|S(\frac\pi k)\phi-\phi\|_1\|x\| + \frac12\omega_x(\frac\pi {|k|}).
\]
Direct computation (as well as \cite[Theorem 3.7]{BK05}) implies that $$\|S(\frac\pi k)\phi-\phi\|_1 \le Const\cdot \frac1{|k|}.$$ Finally, using \eqref{modest} and \eqref{modestt} we get
\eqref{modest1}.
\epf

\bt
Assume that $\X$ is a non-degenerate Banach $L^1(\RR)$-module with the structure associated with a representation $U: \RR \to B(\X)$. Let $B$ be the infinitesimal generator of $U$ and assume that $x\in D(B)$.
Assume also that $y = Bx$ satisfies
\begeq\label{conmodcrit}
\sum_{k\in\NN} \frac{\omega_y(1/k)}{k} <\infty.
\eq
Then $x\in \X_{as}$.
\et

\bpf
For $k\in\NN\backslash\{1\}$, let $f_k\in L^1(\RR)$ be such that $$\hat f_k(\l) = \left\{
\begin{array}{rc}
\frac1{i\l},& \l\ge k-1;\\
\frac{1}{i(2k-2-\l)},& \l < k-1.
\end{array}
\right. 
$$
Similarly, for $-k\in\NN\backslash\{1\}$, let $f_k\in L^1(\RR)$ be such that $$\hat f_k(\l) = \left\{
\begin{array}{rc}
\frac1{i\l},& \l\le k+1;\\
\frac{1}{i(2k+2-\l)},& \l > k+1.
\end{array}
\right. 
$$
Observe that for $k\in\ZZ\backslash\{-1,0,1\}$ we have $\|f_k\|_1 = \frac1{|k|-1}$ and the function
${\varphi}_{k} = f_k*\phi_{k} \in L^1(\RR)$ satisfies 
\[\widehat{\varphi}_{k}(\l) = \frac{\hat\phi_{k}(\l)}{i\l}, \l\in\RR\backslash\{0\}.\]
Hence, for $x\in D(B)$  integration by parts yields $\phi_{k}x = \varphi_{k}(Bx)$, $k\in\ZZ\backslash\{-1,0,1\}$. Finally, using \eqref{modest1}
for $y= Bx$, we get
\[
\|\phi_{k}x\| \le \|f_k\|_1\|\|\phi_{k}y\| \le \frac{Const}{|k|-1}\cdot\omega_y(\frac1{|k|}).
\]
and the result follows.
\epf

Observe that condition \eqref{conmodcrit} is satisfied automatically if $x\in D(B^{1+\a})$, $\a>0$. We also have the following corollary.

\bc
Let $\X = C_{b}(\RR, X)$. Then for any $\a> 0$ we have $C^{1,\a}\subset \X_{as}$.
\ec

\bpf
In this case $B = d/dt$. Since, $x\in C^{1,\a}$ implies \[\omega_{x^\prime}(t) \le \|x\|_{C^{1,\a}}\cdot t^\a,\ t> 0,\] 
the series \eqref{conmodcrit} converges.
\epf

\brem
If $x \in\F=C_0(\RR, X)$ and $x = \hat y$ for some $y\in L^1(\RR)$, then $x\in\F_{as}$ and $\|x\|_{as}\le 5\|y\|_1$.
\erem

We conclude the section with the following useful result.

\bl
If $\F$ is a (spectrally) homogeneous space then $\F_{as}$ is also a (spectrally) homogeneous space.
\el

\bpf
Since convolution operators commute with translation, the verification of the properties of a (spectrally) homogeneous space
is straightforward and is left to the reader.
\epf

\section{Basic properties of the operator $\LL$}\label{operat}

In this section we collect and enhance some of the known spectral properties of  abstract parabolic operators.

Recall that by $A: D(A)\subseteq X\to X$ we denote the infinitesimal generator of a $C_0$-Semigroup $T$.
In a homogeneous space $\F(\RR, X)$ we define the differential operator $\LL = \LL_\F= -d/dt+A$ in \eqref{difop}  as follows
\cite{B96, B13, LZ82}.

\bd\label{defdifop}
A function $x\in\F$ belongs to the domain $D(\LL)$ of the operator $\LL$ if there is a function $y \in\F$ such that
for all $s\le t$ in $\RR$ we have
\begeq\label{semdifop}
x(t) = T(t-s)x(s) - \int_s^t T(t-\tau)y(\tau) d\tau.
\eq
For $x\in D(\LL)$ we let $\LL x = y$, if $x$ and $y$ satisfy \eqref{semdifop}.
\ed
We remark that the operator $\LL$ is well defined as it is not hard to see that for $x\in D(\LL)$ there is a unique
$y$ such that \eqref{semdifop} is satisfied. We also remark that $\F\subset\SS^1$ implies $D(\LL) \subset C_{ub}$.

The operator $\LL$ is \emph{invertible} if it is \emph{injective}, i.e. $\ker\LL = \{0\}$, and \emph{surjective}, i.e.~
its range $\Imm \LL = \LL D(\LL)$ satisfies $\Imm \LL= \F$.

We begin studying the spectral properties of the operator $\LL_\F$ with the following key property of its kernel.
It was originally proved in \cite{B78} for $\F= C_b$.

\bl\label{Linj}
Assume $x\in \ker \LL_\F$. Then
\begeq
i\L(x) \subseteq \s(A)\cap (i\RR).
\eq 
\el

\bpf
Let $x\in \ker \LL_\F$. Then from \eqref{semdifop} we get 
$x(t) = T(t-s)x(s)$, $s\le t$.
Then in view of this equality, for any $f\in L^1(\RR)$, we have
\[
\begin{split}
T(t-s)(f*x)(s) &=  \int_\RR f(s-\tau)T(t-s)x(\tau)d\tau \\ &= \int_\RR f(s-\tau) x(t-s+\tau)d\tau =(f*x)(t).
\end{split}
\]
Hence, \eqref{semdifop} implies $f*x\in \ker\LL_\F$.

Let $i\l_0\notin \s(A)$ and $f\in L^1(\RR)$ be such that $\hat f(\l_0)\neq 0$,  $\supp \hat f$ is compact, and 
$i\supp \hat f \subset \rho(A)$
Then,
according to the above, $y = f*x$ satisfies $y(s+t) = T(t)y(s)$ for all $s\in\RR$ and $t\ge 0$. Observe that since
$f\in C^\infty$, i.e. $f$ is differentiable infinitely many times, we have $y\in C^\infty$. Moreover,
\[y^\prime(s+t) = T(t)Ay(s) = AT(t)y(s),\ s\in\RR, t\ge 0,\]
and, therefore, plugging in $t=0$, we get $y^\prime -Ay = 0$.

Let  $\psi\in L^1(\RR)$ be such that $\hat\psi =1$ in a neighborhood of $\supp \hat f$, $\supp \hat\psi$ is compact, and 
$i\supp \hat\psi \subset \rho(A)$. Consider %the function
$F\in L^1(\RR, B(\F))$ defined by
\[F(t) = \frac1{2\pi}\int_\RR \hat \psi(\l) R(i\l, A)e^{i\l t}d\l, t\in\RR;\]
the above integral makes sense because $i{\supp \hat \psi} \subset \rho(A)$. Observe that 
\begeq
\hat F(\l) = \left\{
\begin{array}{rc}
\hat \psi(\l) R(i\l, A),& i\l\in\rho(A);\\
0, & i\l\notin\rho(A);
\end{array}
\right.
\eq 
is compactly supported and infinitely differentiable. Then, since the operator $A$ is closed,
\[
0 = F*(Ay - y^\prime) = A(F*y)-F^\prime*y -  =  F_0 *y,
\]
where $ F_0 = AF - F^\prime$ has Fourier transform 
$$\widehat { F}_0(\l) = \hat\psi(\l)AR(i\l,A)  - i\l\hat\psi(\l)R(i\l,A) = \hat\psi(\l)I.$$
Hence, $f*x = y =  F_0* y = 0$ and, therefore, $\l_0\notin \L(x)$.
\epf

\bc \label{corinj1}
If the generator $A$ satisfies \eqref{intera} then the operator $\LL_\F$ is injective.
\ec

\bpf
Assume  that $\LL_\F x = 0$. Immediately from Lemma \ref{Linj} we get $\L(x) = \emptyset$. Hence, $x = 0$ and, the operator $\LL_\F$ is injective.
\epf

Next we proceed to use the above result to obtain invertibility conditions for $\LL_\F$ in the case when the
semigroup $T$ is \emph{hyperbolic},
%Assume that the semigroup $T$ is  \emph{hyperbolic},
 i.e.~it satisfies 
\begeq\label{semhyp}
\s(T(1))\cap \TT = \emptyset.
\eq
For such semigroups %are called . For hyperbolic semigroups 
we have
 \[
\s(T(1)) = \s_{in}\cup\s_{out},  
\]
where $\s_{in}$ is the spectral component inside the unit disc and $\s_{out} = \s(T(1))\backslash\s_{in}$.
We let $P_{in}$ and $P_{out}$ be the corresponding spectral projections
\[
P_{in} = \frac1{2\pi i}\int_{\TT} (T(1) - \l I)^{-1} d\l, \quad P_{out} = I- P_{in},
\]
and represent the space $X$ as a direct sum
\[
X = X_{in}\oplus X_{out},\quad X_{in} = P_{in}X,\ X_{out} = P_{out}X.
\]
From the definition of the spectral projections it follows that $P_{in}$ and $P_{out}$ commute with the operators
$T(t)$, $t\ge 0$. Therefore, $X_{in}$ and $X_{out}$ are invariant subspaces for these operators and we can consider the (restriction) semigroup $T_{in}$ and the group $T_{out}$ defined by 
\[
T_{in}: \RR_+ \to B(X_{in}),\  T_{in}(t) = T(t)\vert_{X_{in}};
\]
\[
T_{out}: \RR \to B(X_{out}),\  T_{out}(t) =
\left\{
\begin{array}{ll}
 T(t)\vert_{X_{out}}, & t\ge 0; \\
 (T(-t)\vert_{X_{out}})^{-1};& t< 0.
\end{array} 
 \right.
\]

The following theorem is one of the main results of this paper.
As we mentioned in the introduction, its special cases %of the above result for $\F \in \{C_0,\ \SS^p, C_b\}$ 
appear in \cite{B96, B13, BS10, CL99}.

\bt\label{gfthm}
If $T$ is a hyperbolic semigroup,  the operator $\LL_\F$ from Definition \ref{defdifop} is invertible.  The inverse $\LL^{-1}_\F \in B(\F)$ is defined by
\begeq\label{invop}
(\LL^{-1}_\F y(t) = (G*y)(t) = \int_\RR G(t-\tau)y(\tau)d\tau, \ t\in\RR, \ y\in\F,
\eq
where the Green function $G\in L^1_s(\RR, B(\F))$ is given by
\begeq\label{forgreen}
G(t) =
\left\{
\begin{array}{ll}
 -T(t){P_{in}}, & t\ge 0; \\
 T_{out}(t){P_{out}};& t< 0.
\end{array} 
 \right.
\eq
\et
\bpf

It is immediate from Definition \ref{defhom}\eqref{svert} that the right hand side of \eqref{invop} defines a bound operator in $B(\F)$. We need to check that this operator is, indeed, the inverse of $\LL_\F$.
%Tracing the , we see that it remains valid in a general homogeneous function space: 

%Assume now that $\LL_\F x = 0$. %As we remarked earlier $x\in C_{ub}$. 
Since the semigroup $T$  is hyperbolic,  the spectral inclusion theorem \cite[Theorem IV.3.6]{EN00} %\cite[Lemma 3]{B78} 
implies that the operator $A$ satisfies \eqref{intera}. 
%Lemma \ref{Linj} we deduce
%\[
%i\L(x) \subseteq \s(A) \cap (i\RR) = \emptyset,
%\]
%Since, \eqref{intera} 
%where the last equality holds due to
%\begeq\label{intera1}
%\s(A) \cap (i\RR) = \emptyset
%\eq 
%by 
%. Hence, we get $x = 0$ and, 
Therefore, the operator $\LL_\F$ is injective by Corollary \ref{corinj1}.

%To check that the operator $\LL_\F$ is surjective, given $y\in\F$, we will 
It remains to show that given $y\in\F$ and $x = G*y$ we have $\LL_\F x = y$. %by a direct computation similar to the one in \cite{B96}. ,   
We get
\[
\begin{split}
x(t) &-T(t-s)x(s) = (G*y)(t) - T(t-s)(G*y)(s) \\
& =  \int_t^\infty T_{out}(t-\tau)P_{out}y(\tau)d\tau  -\int_{-\infty}^t T(t-\tau)P_{in}y(\tau)d\tau \\
& - \int_s^{\infty} T(t-s)T_{out}(s-\tau)P_{out}y(\tau)d\tau + \int_{-\infty}^s T(t-\tau)P_{out}y(\tau)d\tau \\
& =  -\int_{s}^t T(t-\tau)P_{out}y(\tau)d\tau -\int_{s}^t T(t-\tau)P_{in}y(\tau)d\tau \\
&= -\int_{s}^t T(t-\tau)y(\tau)d\tau, 
\end{split}
\]
and the result follows from \eqref{semdifop}.
\epf

\bc
If $T$ is a hyperbolic semigroup,  then the generator $A$ satisfies \eqref{intera} and \eqref{resbou}.
\ec

\bpf
As we mentioned above \eqref{intera} follows from the spectral mapping theorem. The inequality \eqref{resbou}
follows since the Fourier transform $\widehat G$ of the Green function $G$ satisfies 
$\widehat G(\l) = R(i\l, A)$, see \cite{BS10} for details. 
\epf

\bc
The operator $\LL_\F$ in Definition \ref{defdifop} is closed.
\ec

\bpf
Let $\omega\in\RR$ be such that the semigroup $T_{\omega}$, $T_\omega(t)= T(t)e^{-\omega t}$, satisfies 
$\|T_\omega(t)\| \to 0$ as $t\to\infty$. The value $\omega_0=\min \omega$, where the minimum is taken
over all $\omega$ with the above property is usually called the \emph{growth bound} of the semigroup $T$ \cite[Definition I.5.6]{EN00}.
 Theorem \ref{gfthm} applied to $T_\omega$ implies that the operator $\LL_\F-\omega I$ is invertible and $(\LL_\F-\omega I)^{-1}y = G_\omega *y$, where
\[
G_\omega(t) =
\left\{
\begin{array}{ll}
 -T_\omega(t), & t\ge 0; \\
 0;& t< 0.
\end{array} 
 \right.
\]
Hence, $\rho(\LL_\F) \neq \emptyset$ and the operator $\LL_\F$ is closed.
\epf

%\newpage

%We begin discussing the spectral properties of $\LL$ with the following result about its resolvent. Recall \cite{EN00} that a $C_0$-semigroup $T$ has type $\omega_0$ if $\|T(t)\| \le Me^{\omega_0 t}$ for some $M>0$ and all $t\ge 0$.

%\bl
%Assume that the operator $\LL_\F= -d/dt+A$ is such that the semigroup $T$  generated by $A$ has type $\omega_0$. Then any $\l\in \CC$ with $\Re e\,\l > \omega_0$ is in the resolvent set of $\LL_\F$ and
%\begeq\label{formobr}
%R(\l, \LL_\F)y = (\LL_\F -\l I)^{-1}y = G_\l*y, \ y\in\F,
%\eq
%where the function $G_\l\in L^1_s(\RR, B(\F))$ is given by
%\begeq\label{GrFu96}
%G_\l(t) = \left\{
%\begin{array}{rl}
% T(t)e^{-\l t}, & t\ge 0; \\
%0,& t< 0.
%\end{array} 
 %\right.
%\eq
%\el

%\bpf ???
%For special classes of homogenous spaces, the proof of the lemma appears in \cite{B96}. Here we observe that
%the same argument applies for all homogeneous spaces given by Definition \ref{defhom}.
%\epf

In what follows we shall denote by $\tilde S: L^1(\RR) \to B(\F)$ the algebra homomorphism given by
\[\tilde S(f)x = f *x, \ f\in L^1(\RR), x\in\F.\]
The next result  asserts that in any homogenous space $\F$ the operator $\LL_\F$
commutes with the operators $\tilde S(f)$, $f\in L^1(\RR)$.

\bl\label{opcom}
For all $f\in L^1(\RR)$ and $x\in D(\LL_\F)$ we have
\[\LL_\F\tilde S(f)x = \tilde S(f)\LL_\F x.\]
\el

\bpf
Assume $\l > \omega_0$ where $\omega_0$ is the growth bound of the semigroup $T$. From \eqref{invop}
we deduce that $R(\l,\LL_\F)\tilde S(f) = \tilde S(f)R(\l,\LL_\F)$. Let $x\in D(\LL_\F)$ and $y = (\LL_\F - \l I)x$. Then
\[\tilde S(f)x = \tilde S(f)R(\l,\LL_\F)(\LL_\F-\l I)x = R(\l, \LL_\F)\tilde S(f)y\]
implies $\tilde S(f)x \in D(\LL_\F)$ and
\[(\LL_\F-\l I)\tilde S(f)x = \tilde S(f)y,\]
from where the result immediately follows.
\epf

\bc
Assume $\Delta \subset \RR$ is closed. Then the spectral submodule $\F(\Delta)$ is an invariant subspace
of the operator $\LL_\F$ and the restriction of $\LL_\F$ to $\F(\Delta)$ coincides with the operator $\LL_{\F(\Delta)}$ given by Definition
\ref{defdifop}.
\ec

\bpf
Assume $x\in \F(\Delta)$ and $y = \LL_\F x$. Let $f\in L^1(\RR)$ be such that $\supp \hat f \cap \Delta = \emptyset$. Then $0 = f*x = \LL_\F(f*x) = f*(\LL_\F y)$ and, therefore, $y\in \F(\Delta)$. 
\epf

Next, we use the above commutativity relation to 
extend the result of Lemma \ref{Linj} to the non-homogeneous case.
%prove the following crucial relation between the Beurling spectra of the functions
%$x$ and $\LL x$ and the spectrum of the generator $A$.

\bl
Assume $x\in D(\LL_\F)$ and $y = \LL_\F x$. Then
\[\L(x) \subseteq \L(y) \cup\{\l\in\RR: i\l\in \s(A)\}.\]
\el

\bpf
Assume $\l\notin \Delta_0 = \L(y) \cup\{\l\in\RR: i\l\in \s(A)\}$ and let $f \in L^1(\RR)$ be such that
$\hat f(\l) \neq 0$, $\supp \hat f$ is compact, and $\supp \hat f\cap \Delta_0 = \emptyset$.
From %Lemma \ref{sprop}(iii) we deduce that $\L(fx)\subseteq \supp\hat f \cap\L(x)$. 
Definition \ref{defdifop}, Lemma \ref{opcom}, and Lemma \ref{sprop}(iv) we deduce that $f*x\in D(\LL_\F)$ and
$\LL_\F(f*x) = f*y = 0$. Hence, %$fx \in C_{ub}\cap \ker \LL$ and \cite[Lemma 3]{B78} 
Lemma \ref{Linj} implies
$i\L(f*x) \subseteq \s(A)\cap(i\RR)$. On the other hand, $\L(f*x)\subseteq \supp\hat f \cap\L(x)$ from
Lemma \ref{sprop}(iii).
 Hence, $\L(f*x) = \emptyset$, $f*x = 0$, and $\l\notin\L(x)$.
\epf

The following corollary is immediate in view of Lemma \ref{sprop}(i).

\bc\label{coring}
Assume $\Delta \subset \RR$ is closed and 
\begeq\label{inter1}
i\Delta \cap \s(A) = \emptyset.
\eq 
Then $\ker \LL_{\F(\Delta)} = \{0\}$.
\ec

\bt\label{main3}
Assume $\Delta \subset \RR$ is compact and the generator $A$ satisfies \eqref{inter1}. Then 
the operator $\LL_{\F(\Delta)}$ is invertible.
\et

\bpf
In view of Corollary \ref{coring} we only need to prove that the operator $\LL_{\F(\Delta)}$ is onto.
Let $y\in \F(\Delta)$, i.e. $\L(y)\subseteq \Delta$. Let also $f\in L^1(\RR)$ be such that
$\hat f = 1$ in a neighborhood of $\Delta$, $i\supp \hat f \subset \rho(A)$ is compact, and $\hat f \in C^\infty$, i.e. $\hat f$  is differentiable infinitely many times. Consider the function $F\in L^1(\RR, B(\F))$ defined by
\[F(t) = \frac1{2\pi}\int_\RR \hat f(\l) R(i\l, A)e^{i\l t}d\l, t\in\RR;\]
the above integral makes sense because $i{\supp \hat f} \subset \rho(A)$. Observe that 
\begeq
\hat F(\l) = \left\{
\begin{array}{rc}
\hat f(\l) R(i\l, A),& i\l\in\rho(A);\\
0, & i\l\notin\rho(A);
\end{array}
\right.
\eq 
is compactly supported and infinitely differentiable, and $F*y \in \F(\Delta)$ by Lemma \ref{sprop}(iii). %We will show that $\LL (F*y) = y$. 
Moreover, since the operators $\LL$ and $A$ are closed, we can
write
\[-\frac d{dt}(F*y)(t) = \frac 1{2\pi}\int_\RR\int_\RR i\l \hat F(\l) e^{i\l(t-s)}y(s)d\l ds\]
and
\[A(F*y)(t) =  \frac 1{2\pi}\int_\RR\int_\RR  \hat f(\l)(I+i\l R(i\l, A))e^{i\l(t-s)}y(s)d\l ds. \]
Hence, $F*y\in D(\LL)$ and $\LL (F*y) = f*y = y$, where the last equality follows from Lemma \ref{sprop}(iv). 
\epf

We conclude this paragraph with a result on the spectrum of the operator $\LL_F$ in a spectrally homogeneous space. For several specific homogeneous spaces this result was proved in \cite{B96}.

\bt
Assume $\F$ is a spectrally homogeneous space. Then
\[\s(\LL_\F) = \s(\LL_\F)+i\RR.\]
\et

\bpf
Directly from Definition \ref{defdifop} and \eqref{repV} we have
\[V(\l)\LL_\F V(-\l) = \LL_\F + i\l I, \ \l \in\RR,\]
and the result follows.
\epf

\section{Invertibility of the operator $\LL$ in $\F_{as}$}\label{sec:4}

The other main result of this paper is the following theorem.

\bt\label{mainthm}
Let $\F = \F(\RR, X)$ be a homogeneous function space and $\F_{as} \subset \F$ be the space of functions with the absolutely summable spectrum. Assume that an operator $\LL = \LL_{\F_{as}}: D(\LL)\subseteq \F_{as}\to \F_{as}$ from Definition \ref{defdifop} is such that
the generator $A$ of the semigroup $T$ satisfies \eqref{intera} and \eqref{resbou}. Then the operator $\LL$ is invertible, $\LL^{-1}\in B(\F_{as})$, and
\begeq\label{invest11}
\|\LL^{-1}\| \le \frac{18}\pi M\left(4+4M+2M^2\right)^{1/2}, \ M   = \sup\limits_{\l\in i\RR} \|R(\l, A)\|.
\eq
\et

The proof of the above result is based on several lemmas.

\bl\label{ml1}
Assume that $\Phi\in C^2(\RR, B(X))$ and $\phi \in L^1(\RR)$ is defined by \eqref{triangle}, i.e.
\[
\hat\phi(\l) = \hat \phi_{0}(\l) = (1-|\l|)\chi_{[-1,1]}(\l),\ \l\in\RR.
\]
Let $\Phi_0 = (\Phi\hat\phi)^{\vee}$, i.e. $\Phi_0$ is the inverse Fourier transform of the function $\Phi\hat\phi$.
Then $\Phi_0\in L^1(\RR, B(X))$ and 
\begeq\label{estfi}
\|\Phi_0\|_1\le \frac2\pi\|\Phi\|_\infty^{1/2}\left(4\|\Phi\|_\infty+4\|\Phi^{\prime}\|_\infty+\|\Phi^{\prime\prime}\|_\infty\right)^{1/2}.
\eq
\el

\bpf
Observe that the definition of  $\Phi_0$ implies that $\Phi_0\in C_b(\RR, B(X))$ and 
\begeq\label{estfi1}
\|\Phi_0(t)\| \le \frac1{2\pi}\|\Phi\|_\infty,\ t\in\RR.
\eq 
We also have
\[\Phi_0(t) = \frac1{2\pi}\int_\RR \Phi(\l)\hat\phi(\l)e^{i\l t}d\l=\frac1{2\pi}\int_{-1}^1 \Phi(\l)(1-|\l|)e^{i\l t}d\l\]
\[
=\frac1{2\pi}\left(\int_{-1}^0 \Phi(\l)(1+\l)e^{i\l t}d\l+\int_0^1 \Phi(\l)(1-\l)e^{i\l t}d\l\right).\]
Applying integration by parts in  the above integrals, we get
\[
\begin{split}
\Phi_0(t) =\frac1{2\pi it}\Big[&\Phi(0) - \int_{-1}^0(\Phi(\l)+\Phi^{\prime}(\l)(1+\l))e^{i\l t}d\l - \\ 
 & \Phi(0) -
\int_{0}^1(-\Phi(\l)+\Phi^{\prime}(\l)(1-\l))e^{i\l t}d\l\Big] \\
 =  \frac1{2\pi it}\Big[&\int_{0}^1\Phi(\l)e^{i\l t}d\l - \int_{-1}^0\Phi(\l)e^{i\l t}d\l - \\
 & \int_{-1}^1 \Phi^\prime(\l)(1-|\l|)e^{i\l t}d\l\Big]. 
\end{split} 
\]
Applying integration by parts once again, we get
\[
\begin{split} 
\Phi_0(t)  &= \frac1{2\pi t^2}\Big[2\Phi(0)-\Phi(1)e^{it}-\Phi(-1)e^{-it}+2\int_{0}^1\Phi^\prime(\l)e^{i\l t}d\l \\
&- 2\int_{-1}^0\Phi^
 \prime(\l)e^{i\l t}d\l-\int_{-1}^1 \Phi^{\prime\prime}(\l)(1-|\l|)e^{i\l t}d\l\Big].
\end{split}
\]
Hence, for $t\neq 0$, we have
\begeq\label{estfi2}
\|\Phi_0(t)\|  \le \frac1{2\pi t^2}\left[4\|\Phi\|_{\infty}+4\|\Phi^\prime\|_{\infty}+\| \Phi^{\prime\prime}\|_{\infty}\right].
\eq
Using \eqref{estfi1} and \eqref{estfi2} we get $\Phi_0\in L^1(\RR, B(X))$ since for any $\a > 0$ 
\[
\begin{split}
\|\Phi_0\|_1 &= \int_{|t|\le \a} \|\Phi_0(t)\| dt + \int_{|t|\ge \a} \|\Phi_0(t)\| dt  \\
&\le \frac{1}\pi\left(\a \|\Phi\|_\infty + \left(4\|\Phi\|_{\infty}+4\|\Phi^\prime\|_{\infty}+\| \Phi^{\prime\prime}\|_{\infty}\right)\int_\a^\infty\frac{dt}{t^2} \right)\\
 &=   \frac{1}\pi\left(\a \|\Phi\|_\infty + \frac1\a\left(4\|\Phi\|_{\infty}+4\|\Phi^\prime\|_{\infty}+\| \Phi^{\prime\prime}\|_{\infty}\right) \right). 
\end{split} 
\]
Plugging in $\a = \left(\frac{4\|\Phi\|_{\infty}+4\|\Phi^\prime\|_{\infty}+\| \Phi^{\prime\prime}\|_{\infty}}{\|\Phi\|_\infty}\right)^{1/2}$ we get \eqref{estfi}.
\epf

\bc
Assume that $\Phi\in C^2(\RR, B(X))$ and $\phi_{n}\in L^1(\RR)$ is defined by \eqref{triangled2} and \eqref{triangle}.
%\[
%\hat\phi_n(\l) = \hat \phi_0(\l+n),\ \l\in\RR,\ n\in\ZZ.
%\]
Let $\Phi_n = (\Phi\hat\phi_{n})^{\vee}$.
Then $\Phi_n\in L^1(\RR, B(X))$ and 
\begeq\label{estfin}
\|\Phi_n\|_1\le \frac2\pi\|\Phi\|_\infty^{1/2}\left(4\|\Phi\|_\infty+4\|\Phi^{\prime}\|_\infty+\|\Phi^{\prime\prime}\|_\infty\right)^{1/2}.
\eq
\ec

\bpf
The result follows by applying the lemma to the function $S(-n)\Phi$, where $S$ is the translation representation
\eqref{transrep}.
\epf

\bc\label{corrn}
Assume that
the generator $A$ of a semigroup $T$ satisfies \eqref{intera} and \eqref{resbou}. 
Let $R_n = (R(i\cdot,A)\hat\phi_{n})^{\vee}$. Then $R_n\in L^1(\RR, B(X))$ and 
\begeq\label{estresn}
\|R_n\|_1\le \frac2\pi M\left(4+4M+2M^2\right)^{1/2}.
\eq
 \ec
 
 \bpf
 The result follows from \eqref{estfin} since $\frac d{d\l} R(i\l, A) = iR^2(i\l,A)$ and $\frac {d^2}{d\l^2} R(i\l, A) = 2R^3(i\l,A)$.
 \epf

 We are now ready to complete the proof of Theorem \ref{mainthm}.

\bpf[Proof of Theorem \ref{mainthm}]
Since $\F_{as}$ is itself a homogeneous space, Corollary \ref{corinj1} applies and the operator $\LL =\LL_{\F_{as}}$ is injective.

To prove surjectivity, consider $y\in \F_{as}$. Let $y_n = \phi_{n}*y$ so that $\L(y_n)\subseteq [n-1, n+1] = \Delta(n)$. Let $R_n$ be defined as in Corollary  \ref{corrn}. Then from the proof of Theorem \ref{main3} and Lemma \ref{sprop}(v) we see that $x_n = R_n * y = R_n *( y_{n-1}+y_n+y_{n+1})\in \F(\Delta(n))$ satisfies $\LL x_n = y_n$. Let $x = \sum_{n\in\ZZ} x_n$, where the series converges absolutely since $y\in\F_{as}$ and $R_n$, $n\in \ZZ$, satisfy \eqref{estresn}. Since the operator $\LL$ is closed, we conclude that $\LL x = y$.

It remains to estimate $\|x\|_{as}$. Observe that 
\[
\phi_{n}*x = \phi_{n}*(x_{n-1} +x_n+x_{n+1}).
\]
Hence, \eqref{estresn} implies
\[
\|\phi_{n}*x\| \le \frac2\pi M\left(4+4M+2M^2\right)^{1/2}\sum_{k = n-1}^{n+1}(\|y\|_{k-1}+\|y\|_k+\|y\|_{k+1}),
\]
and the postulated estimate for $\|\LL\|^{-1}$ follows.
\epf

\brem
Observe that if $\F = C_{ub}$, and $\LL_\F$ satisfies the conditions of Theorem \ref{mainthm}, then $\LL_\F \F \supset \LL_\F \F_{as} \supset  C^{1,\alpha}$, $\a>0$. Moreover, for any function $y\in C^{1,\alpha}$ the equation 
$\LL_\F x = y$ has a unique solution $x\in \F_{as}$.
\erem

\section{Examples}
We begin this section with an example of an operator $A$ and a homogeneous space $\F(\RR, X)$ such that
$A$ satisfies \eqref{intera} and \eqref{resbou} but the operator $\LL=\LL_\F$ is not invertible. This examples 
appears in \cite[Counterexample IV.2.7]{EN00}, we provide it in order to point out a feature that seems to be common for all such examples.

 We let 
$$X = C_0(\RR^+)\cap L^1_\nu(\RR^+),\ \nu(s) = e^s,$$ 
where $L^1_\nu$ is the Beurling algebra of all measurable functions that are summable with the weight $\nu$. The
norm in $L^1_\nu(\RR^+)$ is given by
\[\|x\|_\nu = \int_0^\infty |x(s)|e^sds\]
so that $\|x\|_X = \|x\|_\infty+\|x\|_\nu$.
We let $A = \frac d{ds}$ be the generator of the semigroup $T: \RR_+\to X$ given by $T(t)x(s) = x(s+t)$, $x\in X$,
$s,t\in\RR^+$. Observe that $\|T(t)\| = 1$, $t\ge 0$, the growth bound of $T$ is equal to $0$, and the spectral bound
$s(A) = \sup\{\Re e\l, \l\in\s(A)\}$ satisfies $s(A) \le -1$. Since
\[R(i\l,A)x(s) = \int_0^\infty e^{-\l t}x(s+t)dt,\]
the operator $A$ indeed satisfies 
\eqref{intera} and \eqref{resbou}.
However, if 
we let $\F$ be, for example, $C_0(\RR,X)$, Theorem \ref{thm1} would imply that $\LL_\F$ is not invertible since the spectral radius of the operator $T(1)$ is equal to $1$. Observe that in this case $\s(T(1)) = \TT$. To the best of our knowledge, the semigroups in all of the known examples of this kind have this property.

We conclude the paper with an application of our results to the following  non-linear equation in $\F = C_b$:
\begeq\label{nldif}
x^\prime(t) = (Ax)(t) + y(t) + F(x(t)),\  y\in \F_{as},
\eq
where $F$ is the polynomial
\[F(z) = F_1(z)+F_2(z,z)+\ldots+ F_n(z, z,\ldots, z),\ z\in X,\]
and each $F_k$, $k =1,\ldots, n$, is a $k$-linear map. %By $\|F\|$ we shall mean
%\[
%\|F\| = \|F_1\|+\|F_2\|+\ldots+\|F_n\|.
%\]

We assume that the operator $A$ %generates a hyperbolic semigroup $T$. 
satisfies \eqref{intera} and \eqref{resbou}.
Then $x \in \F_{as}$ is a \emph{mild
solution} of the equation \eqref{nldif} if it satisfies
\[
x = z+\Phi x,
%x(t) = \LL_{\F_{as}}^{-1}y(t)+(\Phi x)(t),
%\int_\RR G(t-s)(y(s) + F(x(s))ds =(G*y)(t) + (G*F(x(\cdot))(t),
\]
%where $G$ is the Green function given by \eqref{forgreen}.
where $z = \LL_{\F_{as}}^{-1}y$ and the non-linear map $\Phi: \F_{as}\to \F_{as}$  is given by 
$\Phi  = \LL_{\F_{as}}^{-1}\circ F$. Observe that the map $\Phi$ is Lipschitz in any ball $B_\beta(0)$ of radius $\beta$ centered
at $0 \in \F_{as}$, that is 
\[
\|\Phi(x) - \Phi(y)\|_{as} \le L_{F,A}(\beta) \|x-y\|_{as},\ x, y\in B_\beta(0).
\]
Moreover, because of \eqref{invest11}, the Lipschitz constant satisfies
\begeq
L_{F,A}(\beta) \le \frac{18}\pi M\left(4+4M+2M^2\right)^{1/2}\sum_{k=1}^n k\|F_k\|\beta^{k-1}. %, \ M   = \sup\limits_{\l\in i\RR} \|R(\l, A)\|.
\eq
\bt
Assume that $\beta > 0$, $y\in\F_{as}$, and $F$ are such that $L_{F,A}(\beta) < 1$ and
$\|\Phi(z)\|_{\F_{as}} < {\beta}(1- L_{A,F}(\beta))$. 
Then the non-linear equation \eqref{nldif} has a unique mild solution 
$x\in \F_{as}$ and $\|x - z\|_{as} \le \beta$.
\et

\bpf 
The solution is obtained by the method of simple iterations as in \cite[Theorem 10.1.2]{D60}.
\epf

\brem
Results of this paper can be extended in a straightforward way to the case of differential inclusions
\[
\frac{dx}{dt}\in\mathscr A x+y,\ x, y\in\F_{as}(\RR, X),
\]
where $\mathscr A$ is a linear relation on $X$ \cite{C98, B08mz, BC02, BK13}.
\erem

%\medskip
%\centerline{\bf Acknowledgements}
%\medskip
%We are grateful to H.~Feichtinger,  K.~Gr\"ochenig and Q.~Sun for inspiration provided by their papers and  useful comments that helped us improve the presentation.

\bibliographystyle{siam}
\bibliography{../refs}

\end {document}